\documentclass[11pt]{amsart} 
\usepackage{amssymb,latexsym, amsmath, amscd, array, graphicx} 

\long\def\forget#1\forgotten{} %

\swapnumbers 
\numberwithin{equation}{section} 

\theoremstyle{plain} 
 
\newtheorem{thm}{Theorem}[section] 
\newtheorem{theorem}[thm]{Theorem}

\newtheorem{prop}[thm]{Proposition} 
\newtheorem{cor}[thm]{Corollary}

\newcommand\theoref{Theorem~\ref} 
 
\newcommand\propref{Proposition~\ref}

\newcommand\remref{Remark~\ref}

\theoremstyle{definition} 
 
\newtheorem{definition}[thm]{Definition}

\newtheorem{remark}[thm]{Remark} 
 
\newtheorem{ex}[thm]{Example}

\newtheorem{conj}[thm]{Conjecture}

\DeclareMathOperator{\inv}{{\rm inv}} 
\DeclareMathOperator{\DR}{{\rm DR}}

\def\1{\hbox{\rm\rlap {1}\hskip.03in{\rom I}}} 
\def\Bbbone{{\rm1\mathchoice{\kern-0.25em}{\kern-0.25em} 
{\kern-0.2em}{\kern-0.2em}I}} 

\long\def\forget#1\forgotten{} % 

\newcommand\ver[1]{\marginpar{\tiny Changed in Ver \VER}} 

\date{\today} 

\begin{document}  

\title [contact and symplectic nilmanifolds ]
{Examples of non-simply connected compact contact manifolds $M$ such that $M \times S^1$ cannot  admit a symplectic structure} 

\author[S. Kutsak]{Sergii Kutsak}

\address {Sergii Kutsak, Department of Mathematics, University of Florida,
1400 Stadium Rd, Gainesville, FL 32611, USA}
\email{sergiikutsak@ufl.edu} 

\begin{abstract}
In this paper we construct non-simply connected contact manifolds $M$ of dimension $\geq5$ such that $M\times S^1$ does not admit a symplectic structure.
\end{abstract}
\maketitle
   
\section{Introduction}

A {\em symplectic structure} on a smooth manifold $ M $ is a non-degenerate skew-symmetric closed 2-form $ \omega \in \Omega^2(M) $. A {\em symplectic manifold} is a pair $ (M,\omega) $ where $ M $ is a smooth manifold and $ \omega $ is a symplectic structure on $ M $. 
The non-degeneracy condition means that for all $ p \in M $ we have the property that if $ \omega(v,w)=0 $ for all $ w \in T_p M $ then $ v=0 $. The skew-symmetry condition means that for all $ p \in M $ we have $ \omega(v,w)=-\omega(w,v) $ for all $ v,w \in T_p M $. The closed condition means that the exterior derivative $ d\omega $ of $ \omega $ is identically zero.  Since each odd-dimensional skew-symmetric matrix is singular, we see that $M$ has even dimension. Every symplectic $2n$-dimensional manifold $ (M,\omega) $ is orientable since the $n$-fold wedge product $ \omega \wedge ...\wedge \omega $ never vanishes. 

A {\em contact structure} on a manifold $M$ of odd dimension $2n+1$ is a completely non-integrable $2n$-dimensional tangent plane distribution $\xi$. In the coorientable case the distribution may be defined by a differential 1-form $\alpha$ as $\xi=\ker \alpha$. Then the non-integrability condition can be expressed by the inequality $ \alpha \wedge (d\alpha)^n \neq 0 $, i.e. the form $ \alpha \wedge (d\alpha)^n$ is nowhere zero. Such differential form $\alpha$ is called a { \em contact differential form}. A contact structure can be viewed as an equivalence class of contact differential forms, where two forms $\alpha$ and $\alpha'$ are equivalent if and only if $\alpha'=f\alpha$ where $f$ is a nowhere zero smooth function on $M$. A { \em contact manifold} is a pair $(M^{2n+1},\eta)$ where $M$ is a smooth manifold of dimension $2n+1$ and $\eta$ is a contact structure on $M$,~\cite{MS}. 
\footnote{MSC (2010): Primary 57R17; Secondary 57R19} 
\footnote{Keywords: symplectic structure, contact structure, nilmanifold, nilpotent Lie algebra, Chevalley-Eilenberg complex} \par
  Every contact manifold $(M,\eta)$ can be embedded as a hypersurface in an exact symplectic manifold. Let $\alpha$ be a contact form on $M$. Then $N=M\times\mathbb{R}$ is a symplectic manifold with symplectic form
 $$ \omega=e^\theta(d\alpha-\alpha\wedge d\theta)=d\lambda, \lambda=e^\theta\alpha$$
 where $\theta$ denotes the coordinate on $\mathbb{R}$, \cite{A}. The manifold $(M\times\mathbb{R}, d(e^\theta\alpha))$ is called a symplectization of $M$. It is natural to ask if there exists a contact manifold $M$ such that $M\times S^1$ cannot admit a symplectic structure? The answer to this question is positive.
 \begin{ex}
 Let $M=S^{2n+1}$ be a sphere of dimension $2n+1$, $n>0$. It is known that $M$ is a contact manifold, \cite{MS}. Clearly, $M \times S^1$ cannot admit a symplectic structure because $H^2(M\times S^1,\mathbb{R})= 0$.
 \end{ex}
 MAIN PROBLEM.   {\em  Find a non-simply connected contact manifold $M$ of dimension $\geq5$  such that $M\times S^1$ does not admit a symplectic structure.}\\
 
 In this paper we give examples of contact manifolds $M$ with nontrivial first homology group $H_1(M)$ such that $M \times S^1$ cannot admit a symplectic structure. In particular, we describe all nilmanifolds $M$ of dimension $5,7$ such that $M \times S^1$ does not admit a symplectic structure. 
 \begin{definition}
A { \em nilmanifold} $M$ is a compact homogeneous space of the form $M=N/\Gamma$, where $N$ is a simply connected nilpotent Lie group and $\Gamma$ is a discrete cocompact subgroup in $N$, ~\cite{TO}. 
\end{definition}
For example, if we consider the group $U_n(\mathbb{R})$ of upper triangular matrices having 1s along the diagonal then the quotient $ M=U_n(\mathbb{R})/U_n(\mathbb{Z}) $ is a nilmanifold, called the Heisenberg nilmanifold, where $U_n(\mathbb{Z}) \subset U_n(\mathbb{R}) $ is the set of matrices having integral entries.\\
Every nilmanifold is the Eilenberg-MacLane space $ K(\Gamma,1) $ for some  finitely presented  nilpotent torsion free group. By Malcev's theorem a discrete group $G$ can be realized as the fundamental group of a nilmanifold if and only if it is a finitely presented  nilpotent torsion free group ~\cite{Ma}. 
\begin{theorem}[\cite{Ma}]\label{t:ma}
A simply connected nilpotent Lie group $N$ admits a discrete cocompact subgroup $\Gamma$ if and only if  there exists a basis $\{e_1, e_2, \dots, e_n\}$ of the Lie algebra $\mathfrak{n} $ of $N$ such that the structural constants $c_{ij}^k$ arising in brackets 
$$ \displaystyle{[e_i, e_j]=\sum_k{c_{ij}^ke_k}} $$
are rational numbers for all $i,j,k$. 
\end{theorem}

  In order to determine if a nilmanifold $N/\Gamma$ admits a symplectic structure  we shall use the fact that the Chevalley-Eilenberg complex of the Lie algebra $\mathfrak{n}$ of a Lie group $N$ is isomorphic to the minimal model of the nilmanifold $N/\Gamma$, \cite{N},\cite{TO}. A nilmanifold $N/\Gamma$ admits a symplectic structure if and only if   the Lie algebra $\mathfrak{n}$ of a Lie group $N$ is symplectic, \cite{KGM}.  A Lie algebra $\mathfrak{n}$ of dimension $2n$ is symplectic if there exists an element $\omega$ of degree $2$ in the Chevalley-Eilenberg complex $(\Lambda \mathfrak{n}^*,\delta)$ of $\mathfrak{n}$ such that $\omega^n \neq 0$ and $\delta\omega=0$.\\
  We give complete lists $\mathfrak{L}_5$, $\mathfrak{L}_7$ of all nilpotent Lie algebras $\mathfrak{g}$ of dimensions $5,7$, respectively, such that $\mathfrak{g}\times\mathfrak{a}$ cannot admit a symplectic structure, where $\mathfrak{a}$ is a one-dimensional Lie algebra. In ~\cite{K} the author provides a complete list $\mathfrak{L}$ of all 7-dimensional nilpotent Lie algebras that admit a contact structure. Note that a nilmanifold $N/\Gamma$ admits an invariant contact structure if and only if $\mathfrak{n}$ admits a contact structure.  Let $M=N/\Gamma$ be a nilmanifold of dimension 7 and let $\mathfrak{n}$ be the Lie algebra of the Lie group $N$. If $\mathfrak{n} \in \mathfrak{L}\backslash\mathfrak{L}_7$ then a nilmanifold $M$ admits an invariant contact structure, $M\times S^1$ cannot admit a symplectic structure and $H_1(M) \neq 0$. Note that there are countably infinitely many nonisomorphic nilpotent Lie algebras which belong to $\mathfrak{L}\backslash\mathfrak{L}_7$.

\section{Acknowledgements}
I would like to thank Alexander Dranishnikov and Yuli Rudyak for many useful discussions.

  \section{Preliminaries}  

   \begin{definition}
 For a Lie algebra $\mathfrak{g}$ the upper central series is the increasing sequence of ideals  defined as follows:
$$ C_0(\mathfrak{g})=0, C_{i+1}(\mathfrak{g})=\{x \in \mathfrak{g}| [x,\mathfrak{g}] \subset C_i(\mathfrak{g})\}. $$
In other words $ C_{i+1}(\mathfrak{g}) $ is the inverse image under the canonical mapping of $ \mathfrak{g} $ onto $ \mathfrak{g}/C_i(\mathfrak{g}) $ of the center of $ \mathfrak{g}/C_i(\mathfrak{g}) $.
 The Lie algebra $ \mathfrak{g} $ is called { \em nilpotent} if there is an integer $k$ such that $ C_k(\mathfrak{g})= \mathfrak{g}.$ The minimal such $k$ is called the { \em index of nilpotency}, ~\cite{B}.
\end{definition}

Let $ \mathfrak{g} $ be a Lie algebra with a basis $ \{X_1, \dots ,X_n\} $. Denote by $ \{x_1,\dots ,x_n\} $ the basis for $ \mathfrak{g}^* $ dual to $ \{X_1, \cdots ,X_n\} $. We obtain a differential $d$ on the exterior algebra $\Lambda\mathfrak{g}^*$
by defining it on degree 1 elements as 
$$ dx_k(X_i,X_j)=-x_k([X_i,X_j]) $$
and extending to $\Lambda\mathfrak{g}^*$ as a graded bilinear derivation. Then
$$ 
\displaystyle{[X_i,X_j]=\sum_{l}{c_{ij}^lX_l}} 
$$
where $c_{ij}^l$ are the structure constants of $ \mathfrak{g} $, and it follows from duality that 
$$ 
dx_k(X_i,X_j)=-c_{ij}^k.
$$
Hence on generators the differential is expressed as
$$ 
dx_k=-\sum_{i<j}{c_{ij}^kx_ix_j}. 
$$
Note that the condition $d^2=0$ is equivalent to the Jacobi identity in the Lie algebra. We call the differential graded algebra
$$ (\Lambda\mathfrak{g}^*,d) $$
the {\em Chevalley-Eilenberg } complex of the Lie algebra $\mathfrak{g} $.\\

\begin{definition}
A differential graded algebra (DGA) $ (\mathcal{M},d) $ is called a {\em minimal model} of a DGA$(A,d_A)$ if the following conditions are satisfied:
\\ (i) $( \mathcal{M},d)=(\Lambda V,d) $ is a free model for $(A,d_A)$;
\\ (ii) $d$ is indecomposable in the following sense: there exists a basis $\{x_\mu: \mu \in I \} $ for some well-ordered index set $ I $ such that $\rm deg (x_{\lambda})\leq \rm deg (x_{\mu}) $ if $ \lambda<\mu $, and each $ d(x_{\lambda}) $ is expressed in terms of $ x_{\mu} $,  $\mu<\lambda$. It follows that $dx_{\mu}$ does not have a linear part.
\end{definition}

\begin{definition}
A minimal model of a smooth manifold $M$ is a minimal model of its de Rham DGA.
\end{definition}

We need the following fact,~\cite{CE}. Let $M=N/\Gamma$ be a nilmanifold. Then the complex of differential forms on $M$ can be identified with the the complex 
of differential forms on $N$ which are right invariant by the elements of $\Gamma$.

\begin{theorem}[\cite{N}]
Let $M=N/\Gamma$ be a nilmanifold.
The natural inclusion of the complex of right invariant differential forms on $N$ into the complex of the differential forms on $N/\Gamma$
$$ \Omega_{\DR}^{\inv}(N) \to \Omega_{\DR}(N/\Gamma)$$
induces an isomorphism on the cohomology level.
\end{theorem}

\begin{cor}[\cite{TO}]\label{c:ce}
The minimal model of a compact nilmanifold $N/\Gamma$ is isomorphic to the Chevalley-Eilenberg complex of the Lie algebra 
$ \mathfrak{n} $ of $ N $.
\end{cor}

 To find all nilpotent Lie algebras $\mathfrak{g}$ of dimension 5,7 such that $\mathfrak{g}\times \mathfrak{a}$ cannot admit a symplectic structure we used classification of indecomposable nilpotent Lie algebras  over $\mathbb{R}$ up to dimension $7$. Note that all nilpotent Lie algebras of dimension up to 6 have been classified, ~\cite{Gr}, ~\cite{BM}. In dimension $5$ there are $9$ nonisomorphic nilpotent Lie algebras and there are $3$ of them that admit a contact structure, ~\cite{K}.

 Many attempts have been done on classification of indecomposable nilpotent Lie algebras of dimension $7$. There are a few lists available: Safiullina  (1964, over $\mathbb{C}$)~\cite{S}, Romdhani (1985, over $\mathbb{R}$ and $\mathbb{C}$)~\cite{R}, Seeley (1988, over $\mathbb{C}$)~\cite{Se1}, Ancochea and Goze (1989, over $\mathbb{C}$)~\cite{AG}. The lists above are obtained by using different invariants. Carles ~\cite{C} introduced a new invariant - the weight system, compared the lists of Safiullina, Romdhani and Seeley, and found mistakes and omissions in all of them. Later in 1993 Seeley incorporated all the previous results and published his list over $\mathbb{C}$,~\cite{Se2}. In 1998 Gong used the Skjelbred-Sund method to classify all 7-dimensional nilpotent Lie algebras over $\mathbb{R}$, ~\cite{G}. We will use Gong's classification with some corrections from the list of Magnin,~\cite{M}. 
Over reals there are 138 non-isomorphic 7-dimensiomal indecomposable nilpotent Lie algebras and, in addition, 9 one-parameter families.
For the one-parameter families, a parameter $\lambda$  is used to denote a structure constant
that may take on arbitrary values (with some exceptions) in $\mathbb{R}.$  An invariant $ I(\lambda) $ is given
for each family in which multiple values of $\lambda$ yield isomorphic Lie algebras, i.e., if $ I(\lambda_1) = I(\lambda_2)$,
then the two corresponding Lie algebras are isomorphic, and conversely. \\

\section{Symplectic structures on the product of 5-dimensional nilpotent Lie algebras with one-dimensional Lie algebra}

Note that there are only finitely many isomorphism classes of nilpotent Lie algebras of dimension less than or equal to $6$. But in higher dimensions there are infinite families of pairwise nonisomorphic nilpotent Lie algebras. In dimension $7$ there are $9$ infinite families which can be parametrized by a single parameter. 
For a Lie algebra $\mathfrak{g}$ of dimension $n$ we denote by $\{x_1,\ldots, x_n \}$ a basis of the dual  space $\mathfrak{g}^*$ and by $\{y\}$ we denote a basis of $\mathfrak{a}^*$ where $\mathfrak{a}$ is a one-dimensional Lie algebra. \par
There are 7 non-isomorphic 5-dimensiomal indecomposable nilpotent real Lie algebras $\mathfrak{g}$ such that $\mathfrak{g}\times\mathfrak{a}$ admits a symplectic structure.

\begin{theorem}\label{t:main1} The following is a complete and non-redundant list $\mathfrak{L}_5$ of all $5$-dimensional nilpotent Lie algebras $\mathfrak{g}$ such that $\mathfrak{g} \times \mathfrak{a}$ admits a symplectic structure:
\vspace{0.2cm}

\rm

\begin{center}
\begin{tabular}{ |l|l|}  
\hline

A$_5$   $ x_1x_2+x_3x_4+x_5y $ & L$_3 \oplus$ L$_2$   $x_2x_5+x_3x_4+x_1y $     \\    \hline  
L$_{5,2}$ $ x_2x_4+x_3x_5+x_1y$ & L$_4 \oplus L_1$  $ x_1x_5+x_2x_4+x_3y $ \\ \hline
L$_{5,3}$  $ x_1x_5+x_3x_4+x_2y $ & L$_{5,4}$  $ x_2x_5-x_3x_4+x_1y $ \\ \hline
L$_{5,6}$  $ x_2x_5-x_3x_4+x_1y $ &   \\ \hline

\end{tabular} 
\end{center}
\end{theorem}

Hence the only 5-dimensional nilpotent Lie algebras $\mathfrak{g}$ such that $\mathfrak{g} \times \mathfrak{a}$ cannot admit a symplectic structure are L$_{5,1}$ and   L$_{5,5}$.  Note that L$_{5,1}$ is the Heisenberg algebra. The nontrivial Lie brackets of the 5-dimensional Heisenberg algebra $\mathfrak{h}$ are defined as follows: \\
$$
 [X_2,X_3]=X_1, [X_4,X_5]=X_1.
$$
Then we can find the differential of the Chevalley-Eilenberg complex of the Heisenberg algebra: 
 $$dx_1=x_2x_3+x_4x_5, dx_i=0, i=2,3,4,5. $$
We want to show that $\mathfrak{h} \times \mathfrak{a}$ cannot admit a symplectic structure. If $\omega$ is a symplectic form on $\mathfrak{h} \times \mathfrak{a}$ then it must contain a term $x_1x_j$ for some $i\in\{2,3,4,5\}$. Assume that $\omega$ contains $x_1x_2$. Note that$d(x_1x_2)=x_2x_4x_5$. Since $d\omega=0$ then $\omega$ must contain a term $x_ix_j$ such that  $d(x_ix_j)$ contains $x_2x_4x_5$. There are no terms $x_ix_j$ such that $d(x_ix_j)$  contains $x_2x_4x_5$  except for  $x_1x_2$. Hence $\omega$ cannot contain a term $x_1x_2$. Similarly one can show that $\omega$ cannot contain $x_1x_3, x_1x_4, x_1x_5$. This is a contradiction. Thus 
$\mathfrak{h} \times \mathfrak{a}$ cannot admit a symplectic structure. It follows that the product of the 5-dimensional Heisenberg manifold $\mathbb{H}_5$ with a circle $S^1$ does not admit a symplectic structure. \par
Note that the product of the 3-dimensional Heisenberg manifold $\mathbb{H}_3$ with a circle $S^1$ is the Kodaira-Thurston manifold $KT$ and it is known that $KT$ is a symplectic manifold, \cite{TO}.
\section{Symplectic structures on the product of 7-dimensional nilpotent Lie algebras with one-dimensional Lie algebra}
 Here we use the same notation for the Lie algebras as in ~\cite{G}, where the Lie algebras are listed in accordance with their upper central series dimensions.  For instance, the Lie algebras that have the upper central series dimensions 1,4,7 are listed as follows: (147A), (147B), (147C), etc. For the sake of brevity we will drop the sign of the wedge product (for instance, $x_ix_j $ means $ x_i \wedge x_j $). \par
There are 63 non-isomorphic 7-dimensiomal indecomposable nilpotent real Lie algebras $\mathfrak{g}$ and, in addition, 2 one-parameter families such that $\mathfrak{g}\times\mathfrak{a}$ admits a symplectic structure.

\begin{theorem}\label{t:main2} The following is a complete and non-redundant list $\mathfrak{L}_7$ of all indecomposable $7$-dimensional nilpotent Lie algebras $\mathfrak{g}$ such that $\mathfrak{g} \times \mathfrak{a}$ admits a symplectic structure:

\rm 

\begin{center}
\begin{tabular}{ |l|l|}  
\hline
(37A)  $ x_1x_7-x_4x_5+x_3x_6+x_2y$ & (37B)  $ x_2x_7-x_4x_6+x_1x_5+x_3y $ \\ \hline
(37C)  $ x_1x_7-x_4x_5+x_3x_6+x_2y $  & (37D) $ x_1x_7-x_4x_5+x_3x_6+x_2y $ \\ \hline
(27B)  $ x_1x_7+x_2x_4-x_3x_6+x_5y $ &  (257A) $ x_1x_6+x_3x_4-x_5x_7+x_2y $ \\ \hline
(257B)  $ x_1x_6+x_2x_7-x_3x_4+x_5y $ & (257C) $ x_1x_6+x_2x_7+x_3x_4+x_5y $  \\ \hline
(257D) $ x_2x_7+x_3x_4+2x_1x_6+x_5y $ & (257F) $ x_1x_3+x_2x_6-x_5x_7+x_4y $  \\ \hline
(257G) $ x_1x_6+x_4x_7+x_2x_3+x_5y $ & (257H) $ x_1x_6+x_3x_4+x_5x_7+x_2y $  \\ \hline 
(257I)  $ x_1x_6+x_2x_7-x_3x_5+x_4y $ & (257J) $ x_1x_6+x_2x_7+x_3x_4-x_3x_5$ \\ &  $+x_4y $  \\ \hline

\end{tabular}
\end{center}

\begin{center}
\begin{tabular}{ |l|l|}  
\hline

(257L) $ x_1x_2+x_3x_6+x_4x_7+x_5y $ & (247A) $ x_1x_7+x_3x_6+x_4x_5-x_2x_7$  \\ &  $+x_5y$  \\ \hline(247C)  $-x_2x_7+x_4x_5+x_3x_6+x_1y $ & (247D)  $ x_1x_7+x_4x_5+2x_3x_6+x_2y$  \\ \hline
(247H)   $x_1x_7+x_4x_5+2x_3x_6+2x_2x_7$   & (247J) $ x_1x_7+2x_2x_6+2x_3x_7-x_4x_5$ 
 \\  $+x_2y$      & $+x_3y $   \\ \hline
(247K)  $x_1x_7+2x_3x_6+x_4x_5+x_2y $ & (247L) $-x_2x_7+x_3x_6+x_4x_5+x_1y$  \\ \hline
(247O) $-x_2x_7+x_3x_6+x_4x_5+x_1y $ & (247Q) $x_1x_7+2x_3x_6+x_4x_5+x_2y $    \\ \hline
(2457A) $x_1x_7+x_2x_6-x_3x_4+x_5y $ & (2457B)  $x_2x_7-x_3x_4+x_5x_6+x_1y $ \\ \hline
(2457C) $x_1x_6+x_2x_7+x_3x_4-x_2x_6$ & (2457D) $x_1x_6+x_2x_7+x_3x_4-x_2x_6$ 
\\ $+x_5y $  & $+x_2x_4+x_5y $ \\ \hline
(2457E) $x_1x_7+x_2x_6+x_3x_4-x_2x_7$ & (2457F)  $x_1x_7+x_2x_6-x_3x_4+x_5y $  
 \\ $+x_5y $ &
\\ \hline
(2457G) $x_1x_7+x_2x_6+x_3x_4-x_3x_5$ & (2457H) $ x_1x_6+x_2x_7-x_3x_4+x_2x_4$ 
\\ $-x_2x_7+x_5y $ & $+x_5y $ \\
 \hline
 
 (2457I )  $x_1x_6+x_2x_7+x_3x_4+x_2x_4$ & (2457J)$x_1x_6+x_2x_7+x_3x_4+x_2x_4$ 
\\ $-x_2x_6+x_5y $ & $-x_2x_6+x_5y$
 \\ \hline
 
(2457K)  $x_1x_6+x_2x_7-x_3x_4+x_2x_4$ & (2357C) $-x_2x_7+2x_3x_6+x_4x_5+x_1y $ 
\\ $+x_5y $ &
\\ \hline 

(2357D) $-x_2x_7+2x_3x_6+x_4x_5+x_1y$ & (147E)  $x_1x_7+\frac{1}{2}x_4x_6-x_2x_5+x_3y, $  
\\  &  when $\lambda=\frac{1}{2}$   \\
& $x_1x_4+x_3x_7+x_5x_6+x_3y, $
\\  &  when $\lambda=-1$ \\ \hline 

(1357G)  $x_1x_7+x_3x_5+x_4x_6+x_2y $ & (1357H) $x_1x_7+x_3x_5-x_3x_6+2x_2x_7$ \\  & $ +2x_4x_5+x_4x_6+x_2y$  \\ \hline 

(1357L)  $x_1x_7+x_2x_5+\frac{1}{2}x_3x_6+x_4y$  & (1357M) $x_1x_6+x_2x_7-x_3x_5+x_4y $, \\ & when $\lambda=2$   \\ \hline
(1357O) $x_1x_7+x_3x_5+x_4x_6+x_2y $ & (1357S) $x_1x_7+x_4x_6+(2-\lambda)x_3x_5$ \\ & $+x_2y, \lambda\neq -1,2 $  \\ \hline
(1357S$_1$)  $x_1x_7+x_2x_6+x_4x_5-x_3x_5$& (13457D)  $x_1x_7-x_2x_5+2x_3x_4 $ 
 \\  $+x_3y $  & $+x_3x_6+x_6y$
\\ \hline
(12457A) $x_1x_7+x_2x_6-x_3x_4+x_4x_5$ & (12457B) $x_1x_7+x_2x_6-x_3x_4$  
 \\ $+x_3x_5+x_5y $ & $+x_4x_5+x_3x_5+x_5y $ \\ \hline
 (12457E)  $x_1x_7+2x_3x_4+x_4x_5$ & (123457A)  $x_1x_7-x_3x_6+x_4x_5 $ 
 \\   $-x_2x_6+x_5y $ &  $+x_2x_7+x_2y$ 
 
   \\ \hline
(123457B)  $x_1x_7-x_3x_6+x_4x_5$ & (123457D) $x_2x_7-x_3x_6+x_4x_5+x_1y$
\\ $+x_2x_7+x_2x_4+x_2y$ & 

 \\ \hline
(123457E) $x_2x_7-x_3x_6+x_4x_5+x_1y$ & (123457I)  $x_2x_7+x_3x_6-x_4x_5+x_1y $ \\  & when $\lambda=-1$ 
\\ \hline
(1357A) $x_1x_7-x_3x_5+x_4x_6+x_2y$  & (1357D) $x_2x_7-x_3x_6-2x_4x_5+x_1y$ \\ \hline

\hline
(37B$_1$)  $x_1x_7+x_2x_6+x_3x_5+x_4y$ & (37D$_1$)  $x_1x_7+x_2x_6+2x_3x_5+x_4y $ \\ \hline
(257J$_1$)  $x_1x_7+x_2x_6-x_3x_4+x_5y $ & (247H$_1$) $ x_1x_7-2x_2x_7-2x_3x_6$
\\  & $-x_4x_5+x_2y $
\\ \hline
(2357D$_1$)  $x_2x_7+2x_3x_6-x_4x_5+x_1y $ & (1357P$_1$)  $x_1x_7-x_3x_6-x_4x_5+x_2y $ \\ \hline
 \end{tabular}
\end{center} 
 
 \begin{center}
\begin{tabular}{ |l|l|}  
\hline

(1357Q$_1$)  $-x_1x_7+x_2x_7+x_3x_6 $ & (1357QRS$_1$)  $ x_1x_7+\lambda +x_3x_6 $ 
\\ $ -x_3x_5 +x_4x_6+x_4y $ & $+(2\lambda-1) x_4x_5+x_2y $ \\
 & $\lambda\neq 0,1,
 \frac{1}{2}$
  \\ \hline

\end{tabular} 
\end{center}
\end{theorem}

It is straightforward to verify that all 2-forms listed in \theoref{t:main2} are symplectic. We want to show that the Lie algebra (13457C) cannot admit a symplectic structure. We can find the differential of the Chevalley-Eilenberg complex of the Lie algebra (13457C):

\begin{tabular}{ l l l  }  
 $dx_7=x_1x_6+ x_2x_5-x_3x_4$ & $dx_5=x_1x_4,$ & $ dx_4=x_1x_3$ \\
    $dx_3=x_1x_2,$ & $dx_i=0, i=1,2.$ \\
\end{tabular} \\
Assume that the Lie algebra (13457C) admits a symplectic form $\omega$. Then $\omega$ must contain $x_{i_1}x_{i_2}+x_{i_3}x_{i_4}+x_{i_5}x_{i_6}+x_{i_7}y$, where $\{x_{i_1}, \ldots, x_{i_7}\}=\{x_1,\ldots x_7\}$.  In particular, $\omega$ must contain a term of the form $x_jx_7$ for some $j \in \{1,\ldots,7\}$. Suppose that $\omega$ contains a term $x_1x_7$. Note that $d(x_1x_7)=-x_1x_2x_5+x_1x_3x_4$. Since $d\omega=0$ then it must contain terms $x_ix_j$ such that $d(x_ix_j)$ contains $x_1x_2x_5, x_1x_3x_4$. The only such term is $x_3x_5$ and we compute $d(x_3x_5)=x_1x_2x_5+x_1x_3x_4$. Note that $d(x_1x_7+x_3x_5)=2x_1x_3x_4$ and $d(x_3x_5-x_1x_7)=2x_1x_2x_5$. There are no terms $x_ix_j$ such that $d(x_ix_j)$  contains $x_1x_3x_4$ or $x_1x_2x_5$ except for  $x_3x_5, x_1x_7$. Hence a 2-form $\omega$ cannot contain a term $x_1x_7$. Suppose that $\omega$ contains a term $x_2x_7$. Then $d(x_2x_7)=x_1x_2x_6+x_2x_3x_4$ and $d(x_ix_j)$ does not contain $x_2x_3x_4$ for all $i,j$ except for $x_2x_7$. Hence if $\omega$ contains $x_2x_7$ then $d\omega \neq 0$. Similarly it can be shown that $\omega$ cannot contain $x_3x_7, x_4x_7, x_5x_7, x_6x_7$. This is a contradiction. Thus the Lie algebra (13457C) cannot admit a symplectic structure. Similarly one can show that all 7-dimensional nilpotent Lie algebras which are not listed in \theoref{t:main2} cannot admit a symplectic structure. \\

\begin{remark}
In ~\cite{K} the author gives a complete and non-redundant list $\mathcal{L}$ of all indecomposable $ 7$-dimensional nilpotent real Lie algebras that admit a contact structure. In particular, it was proved that all 9 one-parameter families admit a contact structure. Hence there are 7 one-parameter families
of indecomposable $7$-dimensional nilpotent real Lie algebras $\mathfrak{g}$ such that $\mathfrak{g}\times\mathfrak{a}$ does not admit a symplectic structure. Thus we can produce infinitely many nondiffeomorphic nilmanifolds $M=N/\Gamma$ of dimension $7$ such that $M$ admits an invariant contact structure and $M\times S^1$ does not admit a symplectic structure.
\end{remark}

\begin{remark}
In 1992 Seeley has solved the problem of estimating the number of parameters $F_n$ needed to classify the laws of $n$-dimensional complex nilpotent Lie algebras, and came up with the estimation $F_{n+2}\geq n(n-1)(n+4)\backslash 6 - 3$, \cite{Se3}. Hence it is very difficult to write a complete list for dimensions greater than $7$. For example, for dimensions $8$ and $10$ the number of parameters involved will be respectively $\geq4$ and $\geq13$. 
\end{remark}

 \begin{remark}\label{rem:form}
   Note that all symplectic structures listed above are of the form
   $$\displaystyle{\sum_k{x_{i_k}x_{j_k}}+x_ly, \text{ where } dx_l=dy=0}.$$
Hence we can construct a symplectic structure on $\mathfrak{g} \times \mathfrak{h}$ for all Lie algebras $\mathfrak{g}$, $\mathfrak{h}$ in $\mathfrak{L}_5 \cup \mathfrak{L}_7$.   Because if 
$\displaystyle{\sum_k{x_{i_k}x_{j_k}}+x_ly}$ is a symplectic structure on $\mathfrak{g}\times\mathfrak{a}$ and $\displaystyle{\sum_n{y_{s_n}y_{t_n}}+y_ry}$  is  a symplectic structure on $\mathfrak{h}\times\mathfrak{a}$ then
$$\displaystyle{\sum_k{x_{i_k}x_{j_k}}+x_ly_r+\sum_n{y_{s_n}y_{t_n}}}$$
is a symplectic form on $\mathfrak{g} \times \mathfrak{h}$. 
\end{remark}

 For example, let $\mathfrak{g}$ be the Lie algebra (2457G) with basis $\{x_1, \ldots x_7\}$ and let $\mathfrak{h}$ be the Lie algebra (1357Q$_1$) with the basis $\{y_1, \ldots y_7\}$. Then a symplectic structure on $\mathfrak{g} \times \mathfrak{h}$ is given by a form: 
$$\omega=x_1x_7+x_2x_6+x_3x_4-x_3x_5-x_2x_7+x_5y_4-y_1y_7+y_2y_7+y_3y_6-y_3y_5+y_4y_6.$$
Moreover, in a similar fashion we can construct a symplectic structure on $ \mathfrak{g}_1 \times \ldots \times \mathfrak{g}_{2n}$, for all $\{\mathfrak{g}_1, \ldots, \mathfrak{g}_{2n}\} \subset \mathfrak{L}_5 \cup \mathfrak{L}_ 7 $.  One can see that if $\mathfrak{g}\times\mathfrak{h}$ admits a symplectic structure $\omega$ then $\omega$ does not contain a term $x_iy_j$ unless $dx_i=dy_j=0$.

   \begin{prop}\label{p:product3}     For all nilpotent Lie algebras $\mathfrak{g}, \mathfrak{h}$ of dimension 5,7 the following holds true: \\
 (i) $\mathfrak{g} \times \mathfrak{a}$ admits a symplectic structure if and only if $\mathfrak{g} \times \mathfrak{g}$ does so. \\
 (ii) $\mathfrak{g} \times \mathfrak{h}$ admits a symplectic structure if and only if both $\mathfrak{g} \times \mathfrak{g}$ and $\mathfrak{h} \times \mathfrak{h}$  admit a symplectic structure. 
 \end{prop}
 
  \begin{prop}\label{p:product1}
 Let $X$ be a nilmanifold of dimension 5,7. Then the following holds true:\\
 (i) $X\times S^1$ admits a symplectic structure if and only if $X\times X$ does so.\\
 (ii) $X\times Y$ admits a symplectic structure if and only if $X\times X$ and $Y\times Y$ admit a symplectic structure.
\end{prop}

Note that  \propref{p:product3} follows from \remref{rem:form}.  \propref{p:product1} follows from \propref{p:product3}.\\

CONJECTURE. {\em For all odd-dimensional nilpotent real Lie algebras $\mathfrak{g}, \mathfrak{h}$ the following holds true: \\
 (i) $\mathfrak{g} \times \mathfrak{a}$ admits a symplectic structure if and only if $\mathfrak{g} \times \mathfrak{g}$ does so. \\
 (ii) $\mathfrak{g} \times \mathfrak{h}$ admits a symplectic structure if and only if both $\mathfrak{g} \times \mathfrak{g}$ and $\mathfrak{h} \times \mathfrak{h}$  admit a symplectic structure. }

\begin{ex}
Recall that Heisenberg algebra of dimension $2n+1$ is a Lie algebra $\mathfrak{g}$  with Lie brackets $[X_{2i},X_{2i+1}]=X_1$, where $\{X_1,\ldots,X_{2n+1}\}$ is a basis of $\mathfrak{g}$ and $n \geq 2$. We consider Chevalley-Eilenberg complex of $\mathfrak{g}$ and it is not difficult to see that $\eta=x_1$ is a contact form on $\mathfrak{g}$ and that $\mathfrak{g}\times\mathfrak{a}$ cannot admit a symplectic structure. Hence the Heisenberg manifold  $\mathbb{H}_{2n+1}$  of dimension $2n+1 \geq 5$ is a nilmanifold with nontrivial first homology group ${H}_1(  \mathbb{H}_{2n+1})$ which admits an invariant contact structure and such that $\mathbb{H}_{2n+1}\times S^1$  cannot admit a symplectic structure.
\end{ex}

\end{document}